\documentclass[11pt]{article}
\usepackage[english]{babel}
\usepackage{amssymb,amsmath,amsthm}
\usepackage[dvipsnames]{pstricks}
\usepackage{MnSymbol}

\newtheorem{theorem}{Theorem}[section]
\newtheorem{lemma}[theorem]{Lemma}
\newtheorem{proposition}[theorem]{Proposition}

{\theoremstyle{definition}}
{\theoremstyle{definition}}
{\theoremstyle{definition}\newtheorem{remark}[theorem]{Remark}}

\numberwithin{equation}{section}

\def\Z{{\mathbb Z}}

\def\K{{\mathbb K}}

\def\epsilon{\varepsilon}
\def\kappa{\varkappa}

\def\phi{\varphi}
\def\leq{\leqslant}
\def\geq{\geqslant}

\def\dim{{\rm dim}\,}

\title{ On Koszulity in homology of moduli spaces of stable n-pointed curves of genus zero}

\author{Natalia Iyudu}

\date{}

\begin{document}

\maketitle


{\bf Abstract}

We prove Koszulity of the homology of the of moduli spaces of stable n-pointed curves of genus zero
 $\overline{M}_{0,n}$ for $n=5$, using its presentation due to Keel and the Priddy criterion of Koszulity. For $n=6 $ we establish that $\overline{M}_{0,6}^!$ is potential, find the expression for the potential, and based on that prove that it is Koszul.

\section{Introduction}

In the paper \cite{Ke}, due to Keel, dedicated to the intersection theory  of the moduli space of $n$-pointed stable curves of genus zero, it was shown that the canonical map from the Chow ring to the homology is an isomorphism, and the presentation of this ring by generators and relations have been derived.

The variables $D^S$ in this presentation are parameterized by subsets $S \subset \{1,\cdots,n\}$ consisting of two or more elements,  variables corresponding to complimentary sets are coincide: $D^S=D^{S^C}$. The relations then looks as follows:

(1). For any four distinct elements $i,j,k,l \in \{1,\cdots,n\}$:

$$ \sum_{i,j \in S, k,l \notin S} D^S = \sum_{i,k \in S, j,l \notin S} D^S =
\sum_{i,l \in S, j,k \notin S} D^S $$

(2). $ D^S D^T$ unless one of the following holds:

$$ S \subset T, T \subset S, S \subset T^C, T^C \subset S. $$

We will use this presentation to examine Koszulity of the ring $A=H_{\bullet}(\overline{M}_{0,n})$.

For $n=5$ we suggest the following choice  of linearly independent variables: out of all generators

$$\{  \delta_{i,j}   \}=\{  \delta_{12}, \delta_{13},\delta_{14},\delta_{15},\delta_{23},\delta_{24},\delta_{25},\delta_{34},\delta_{35},\delta_{45}   \}$$

we take those, sitting on the sides of the pentagon:

$$  x_i =\delta_{i,i+1}, \,\, i \in \Z / 5\Z,$$

then we rewrite all other variables and relations via them. This presentation wakes up some associations with cluster  algebras. Using the obtained relations we compute a presentation for the Koszul dual (noncommutative) algebra $A^!$. It turns out that it can be presented by one relation of the form:

 $$ \sum_{a\neq b, a \neq b+1, a \neq b-1} x_ax_b - \sum_{k=1,...,5} x_k^2=0.$$

We can find a change of variables, such that the new relations does not have a square term $x^2$. Thus
we have an algebra $A^!$ presented by a quadratic Gr\"obner basis.
Hence due to the Priddy criterion \cite{Pr, Manin}, we can conclude that $A^!$ is Koszul, and hence
$A$ itself as well.

Then we consider $A=H_{\bullet}(\overline{M}_{0,6})$. After change of generators to a smaller set of variables, which span the same space modulo linear relations, we are able to show that the dual algebra $A^!$ is potential, and find an expression for the potential. Then we construct a potential complex, which is coincide  with the Koszul complex. We show its exactness by checking it is exact in all but one places, then numeric Koszulity implies Koszulity, according to lemma\ref{NKos}.


Our methods here could be compared to those appeared in a vast literature, where Koszulity of various algebras given by relation, originated from geometric objects, have been proved. To name a few, there are series of papers \cite{R1, R2, R3} dedicated to Koszulity of algebras associated to graphs and sell complexes, papers \cite{Cox} on Koszulity of objects associated to the Coexter group action, etc.

\section{Koszulity of $H_{\bullet}(\overline{M}_{0,5})$}

We consider first the presentation of $A=H_{\bullet}( \overline{M}_{0,n})$ due to Keel \cite{Ke} with generators

$$\{  \delta_{i,j}   \}=\{  \delta_{12}, \delta_{13},\delta_{14},\delta_{15},\delta_{23},\delta_{24},\delta_{25},\delta_{34},\delta_{35},\delta_{45}   \}$$

and corresponding linear and monomial quadratic relations. Let $V={\rm span}_k  \{\delta_{i,j}\}$. We suggest to choose
as a linear basis in $V$ the following set: $\{  x_i =\delta_{i,i+1}, i \in \Z / 5\Z   \}$. Via this basis we are going to express other elements of $\{  \delta_{i,j}   \}$, then we rewrite quadratic relations in terms of variables $x_i, i \in \Z / 5 \Z $. This will allow us later on to pass to the Koszul dual algebra $A^!$, which could be shown to be Koszul. The presentation for the $A^!$ which we obtain, after suitable change of variables, will have a quadratic Gr\"obner bases. This due to the  Priddy's PBW criterion \cite{Pr, Manin}  will imply Koszulity of $A^!$ and hence of $A$.

\vspace{3mm}

{STEP 1}

Expression of $\{  \delta_{i,j}   \}$ via $x_i=\delta_{i,i+1}, i \in \Z / 5 \Z $.

\vspace{3mm}

Linear relation which we have for quadruple $ \{a,b,c,d\} \in \Z / 5 \Z $ of pairwise distinct elements from 1 to 5, are

$$ \delta_{ab}+\delta_{cd}=\delta_{ac}+\delta_{bd}=\delta_{ad}+\delta_{bc}.$$

In particular, $ \delta_{ad}=\delta_{ab}+\delta_{cd}-\delta_{bc}$, and for, say $a=i, b=i+1, c=i+2, d=i+3$
we get

 $$ \delta_{i,i+3}=\delta_{i,i+1}+\delta_{i+2,i+3}-\delta_{i+1,i+2},$$

 that is

$$ \delta_{i,i+3}=x_{i}+x_{i+2}-x_{i+1}.$$

The set of elements $\{  \delta_{i,j}   \}$, which are our original generators can be presented 'geometrically' as a sides and diagonals of the pentagon (analogously to the presentation of cluster variables). Sides correspond to variables $\delta_{i,i+1}$ and diagonals to $\delta_{i,i+3}$, for appropriate $i$. We will call these two types of generators {\it side} type and {\it diagonal} type generators respectively. Hence, by expressing the diagonal type variables $\delta_{i,i+3}$,  via those $\delta_{i,j}$ with nearby indexes, e.i. via the side type variables, we express all generators via the new basis $\{x_{i}, i \in \Z / 5 \Z \}$.

\vspace{3mm}

{STEP 2}

Calculation of quadratic relations in terms of $x_i=\delta_{i,i+1}, i \in \Z / 5 \Z $.

\vspace{3mm}

Our initial quadratic relations have a shape

$$  \delta_{i,j}  \delta_{j,k}, i \neq k  $$

This means, for those variables which are side type, we will have

$$ (1)    \quad  \delta_{i,i+1}  \delta_{i+1,i+2} = x_i x_{i+1} = 0.  $$

If a relation contains both side and diagonal type generators, it have a shape:

$$  \delta_{i,i+3} x_{k},  k \neq i+1, k \in \Z / 5 \Z .  $$

In the pentagon picture it could be interpreted as a fact, that only the side $x_{i+1}$ does not intersect the diagonal
$ \delta_{i,i+3}$.

So in our new variables these relations  look like:


$$ (2)  \quad  (x_i +x_{i+2} - x_{i+1}) x_i = 0 (k=i)$$

$$    (x_i +x_{i+2} - x_{i+1}) x_{i-1} = 0 (k=i-1)$$

$$     (x_i +x_{i+2} - x_{i+1}) x_{i+2} = 0 (k=i+2)$$

$$     (x_i +x_{i+2} - x_{i+1}) x_{i+3} = 0 (k=i+3).$$

\vspace{3mm}

In the case when relations contain two diagonal variables, we also will have:

 $$\delta_{i,i+3}\delta_{j,j+3}=0, $$ which means

  $$\delta_{i,i+3}\delta_{i-3,i}=0.$$

 That is,

 $$  (3)   \quad  (x_i +x_{i+2} - x_{i+1}) (x_{i-3} +x_{i-1} - x_{i-2}) = 0.$$

 Since it is obviously the same as

 $$    (x_i +x_{i+2} - x_{i+1}) (x_{i+2} +x_{i-1} - x_{i+3}) = 0,$$

 we can see that (3) follows from (2) and hence the relations of $A$ on variables $x_i$ are formed by (1) and (2).

 This system of relations is equivalent to the following one:

 $$ x_ix_{i+1}=x_{i+1}x_i=0$$

 and for any other pair $a,b,a',b' \in \Z / 5 \Z $

 $$ x_ax_b=x_{a'}x_{b'}=-x_a^2.$$

 The Hilbert series here is clearly $H_A=1+5t+t^2$. This also implies  that codimension of the second graded component of the ideal $I_2$ is 1. Now using this system of relations, and additional information that the
 dimension of the space orthogonal to $I_2$ is one, we can write down generating relation for the Koszul dual algebra $A^!$:

 $$ \sum_{a\neq b, a \neq b+1, a \neq b-1} x_ax_b - \sum_{k=1,...,5} x_k^2=0.$$

 This generating system of $A^!$ is not yet form a Gr\"obner basis of the ideal. But we can make the change of variables
 which will make this relation free from square of one of variables (say, $x_2$).
 This will make a presentation of algebra combinatorially free, that is new relations is a Gr\"obner basis. So the Gr\"obner basis is quadratic
 Due to the Priddy PBW criterion \cite{Pr} this implies that $A^!$ is Koszul.

So we have proved

\begin{theorem} The algebra $A=H_{\bullet}(\overline{M}_{0,5})$  is Koszul.
\end{theorem}

\vspace{3mm}

\section{Koszulity in the case of $\overline{M}_{0,6}$}

Our goal in this section will be the following.

\begin{theorem}\label{theorem00} The algebra $A=H_{\bullet}(\overline{M}_{0,6})$  is Koszul.
\end{theorem}

The canonical generating set for $A$ from the described above Keel's presentation we denote $\{a_{ij}=a_{\{i,j\}},S_{ijk}=S_{\{i,j,k\}}\}$ for two-element subsets $\{i,j\}$ and three-element subsets $\{i,j,k\}$ of $\Omega=\{1,2,3,4,5,6\}$.. We have $15$ generators $a_{ij}=a_{ji}$ with $i\neq j$ and 10 generators $S_{R}=S_{\Omega\setminus R}$ for $3$-element subsets $R$ of $\Omega$.

It will be convenient to use the following notation:
$$
\begin{array}{l}
ijR\iff i,j\in R\ \text{or}\ i,j\in\Omega\setminus R;
\\
iRj\iff i\in R,\ j\in\Omega\setminus R\ \text{or}\ j\in R,\ i\in\Omega\setminus R.
\end{array}
$$
The defining relations of $A$ are as follows.

The linear relations are:
$$
0=\eta_{i,j,k,m}=a_{ij}+a_{km}+S_{ij\alpha}+S_{ij\beta}-a_{ik}-a_{jm}-S_{ik\alpha}-S_{ik\beta},
$$
for all $i,j,k,m,\alpha,\beta\in \Omega$ such that $\{i,j,k,m,\alpha,\beta\}=\Omega$.

The quadratic relations are:
$$
\begin{array}{l}
S_RS_{R'}=0\ \ \text{if}\ \ S_R\neq S_{R'};
\\
S_Ra_{ij}=0\ \ \text{if}\ \ iRj,
\\
a_{ij}a_{ik}=0\ \ \text{if}\ \ i,j,k\ \text{are distinct}
\\
\text{+commutativity}.
\end{array}
$$

\begin{remark}\label{rem1} Note that there are automorphisms of $A$ corresponding to elements of the group $S_6$ of permutations of points $\Omega$. Indeed, the natural action of $S_6$ on the set $\{a_{ij},S_{ijk}\}$ of generators extends to action by automorphisms on $A$. This follows from the fact that this action on the set of generators preserves the way in which corresponding sets intersect and therefore preserves the set of defining relations.
\end{remark}

Consider the set of variables
$$
v_i=\sum_{j\neq i} a_{ij}.
$$
Our goal is to show that

\begin{theorem}\label{thm01} The linear span of all $a_{ij}$ and $S_R$ and the linear span of all $v_i$ and $S_R$ coincide modulo $L$, where $L$ is the linear span of all linear relations $\eta_{i,j,k,m}$. Furthermore, $S_R$ and $v_k$ are linearly independent modulo $L$. Finally, the expressions of $a_{ij}$ via $v_i$ and $S_R$ modulo $L$ are:
$$
a_{ij}=\frac15(v_i+v_j)-\frac1{20}\sum_{k\neq i,j}v_k+\frac15\sum_{R:iRj}S_R-\frac3{10}\sum_{R:ijR}S_R.
$$
\end{theorem}

First, we need the following lemma.

\begin{lemma}\label{perp} The space $L^\perp$ is spanned by $v_i$ for $i\in\Omega$ and $\xi_R=S_R-\frac12\sum\limits_{pqR}a_{pq}$. Furthermore, the six vectors $v_i$ and $10$ vectors $\xi_R$ form a linear basis in $L^\perp$.
\end{lemma}

\begin{proof}First, note that introduced above new variables $v_s$ are orthogonal to each $\eta_{i,j,k,m}$. This gives us six linearly independent elements $v_i$ of $L^\perp$. Next, we search for vectors orthogonal to $L$ of the form $\xi_R=S_R+\alpha\sum\limits_{pRq}a_{pq}+\beta\sum\limits_{pqR}a_{pq}$. The scalar products of such a vector with vectors $\eta_{i,j,k,m}$ have to be equal  to 0, we see that $w\in L^\perp$ precisely when $1+2\beta-2\alpha=0$. Now each $\xi_R$ belongs to $L^\perp$. Obviously, all 16 vectors $v_i$ and $\xi_R$ are linearly independent. On the other hand, the Hilbert series of $A$ is known \cite{Manin}: $H_A=1+16t+16t^2+t^3$ and therefore, the dimension of $L$ is 25-16=9 (25 is the number of generators). Hence, the dimension of $L^\perp$ is $25-9=16$. Hence the 16 vectors $v_i$ and $\xi_R$ form a linear basis in $L^\perp$.
\end{proof}

\begin{proof}[Proof of Theorem~\ref{thm01}] We shall express $a_{ij}$ as a linear combination of $v_k$ and $S_R$ modulo $L$. We are going to use the basis in $L^\perp$ from Lemma~\ref{perp}. We shall try to find $\alpha,\beta,\gamma,\delta\in\K$ such that
$$
a_{ij}=\alpha(v_i+v_j)+\beta\sum_{k\neq i,j}v_k+\gamma\sum_{R:iRj}S_R+\delta\sum_{R:ijR}S_R
$$
modulo $L$. This is equivalent to the inclusion
$$
\zeta_{ij}=a_{ij}-\alpha(t_i+t_j)-\beta\sum_{k\neq i,j}v_k-\gamma\sum_{R:iRj}S_R-\delta\sum_{R:ijR}S_R\in L.
$$
This, in turn, is equivalent to $\zeta_{ij}$ being orthogonal to all $v_k$ and $\xi_R$ of Lemma~\ref{perp}. Orthogonality of $\zeta_{ij}$ to $v_i$ (or $v_j$) reads as the equation $1-6\alpha-4\beta=0$, while orthogonality of $\zeta_{ij}$ to $v_k$ with $k\notin\{i,j\}$ is equivalent to the equation $2\alpha+8\beta=0$. Orthogonality of $\zeta_{ij}$ to $\xi_R$ satisfying $ijR$ reads as the equation $-1+4\alpha+8\beta-\delta=0$. Finally, orthogonality of $\zeta_{ij}$ to $\xi_R$ satisfying $iRj$ is equivalent to the equation $4\alpha+12\beta-\gamma=0$. Solving this system of linear equations, we obtain that $\zeta_{ij}\in L$ precisely when $\alpha=\frac15$, $\beta=-\frac1{20}$, $\gamma=\frac15$ and $\delta=-\frac3{10}$. Hence
$$
a_{ij}=\frac15(t_i+t_j)-\frac1{20}\sum_{k\neq i,j}v_k+\frac15\sum_{R:iRj}S_R-\frac3{10}\sum_{R:ijR}S_R
$$
modulo $L$.

All other statements of Theorem~\ref{thm01} follow via elementary dimension arguments..
\end{proof}

The above results immediately yield the following presentation of $A$.

\begin{theorem}\label{thm02} $A$ is given by $16$ generators $v_i$, $S_R$ and quadratic relations$:$
$$
\begin{array}{l}
S_RS_{R'}=0\ \ \text{if}\ \ S_R\neq S_{R'};
\\
S_Rt_{ij}=0\ \ \text{if}\ \ iRj,
\\
t_{ij}t_{ik}=0\ \ \text{if}\ \ i,j,k\ \text{are distinct}
\\
\text{+commutativity},
\end{array}
$$
where
$$
t_{ij}=4(v_i+v_j)-\sum\limits_{k\neq i,j} v_k+4\sum\limits_{R:iRj} S_R-6\sum\limits_{R:ijR} S_R.
$$
\end{theorem}

\begin{theorem}\label{thm03}
The dual algebra $B=A^!$ is potential with symmetric potential which is also invariant under the $S_6$-action corresponding to permutations of $\Omega$. Furthermore, the corresponding potential, when written in terms of generators $v_i$ and $S_R$ has the form

\begin{align*}
&P=a_1\sum_R S_R^3+a_2\sum_{R,i} S_R^2v_i^\rcirclearrowleft+a_3\sum_{R,i}S_R{v_i^2}^\rcirclearrowleft+a_4\sum_{i,j,R:ijR} S_Rv_iv_j^\rcirclearrowleft
\\
&+a_5\sum_{i,j,R:iRj} S_Rv_iv_j^\rcirclearrowleft+a_6\sum_{i} v_i^3+a_7\sum_{i,j:i\neq j} {v_iv_j^2}^\rcirclearrowleft+a_8\sum_{i,j,k\atop i,j,k\ \text{are distinct}} v_iv_jv_k
\end{align*}

for some non-zero $(a_1,\dots,a_8)\in \K^8$.
\end{theorem}

\begin{proof} Let $V=A_1$ be the linear span of the generators $v_i$ and $S_R$ of $A$. Let $R$ be the subspace of $V^2$ of quadratic relations of $A^!$. Then $R^\perp$ is the space of quadratic relations of $A$.
Since
the Hilbert series of $A$ is $H_A=1+16t+16t^2+t^3$ \cite{Manin},, $\dim A_3=1$. Hence the codimension of $VR^\perp+R^\perp V$ in $V^3$ equals $1$. It follows that $(VR^\perp+R^\perp V)^\perp=RV\cap VR$ is one-dimensional. Let $P$ be a (unique up to a non-zero scalar multiple) element of $V^3$ spanning $RV\cap VR$. The fact that $R^\perp$ contains commutators of generators easily implies that $P$ is symmetric. In particular it is cyclicly symmetric, aka a potential. The fact that $P$ spans $RV\cap VR$ together with $H_A=1+16t+16t^2+t^3$ now yields that $R$ is spanned by the first derivatives of $P$. Since $R^\perp$ is invariant under $S_6$-action, so is $R$ and so is $RV\cap VR$. It follows that $P$ is invariant under this action. Finally, since $S_RS_{R'}\in R^\perp$, we have that two distinct $S_R$ never feature in the same monomial of $P$. It follows that $P$ must have the form as stated.
\end{proof}

Now we shall compute $P$.

\begin{theorem}\label{thm04}
The dual algebra $B=A^!$ is potential with potential $($in generators $v_i$, $S_R)$ given by the formula
\begin{align*}
&P=2\sum_R S_R^3-2\sum_{R,i} S_R^2v_i^\rcirclearrowleft+4\sum_{i,j,R:iRj} S_Rv_iv_j^\rcirclearrowleft
\\
&+5\sum_{i} v_i^3-11\sum_{i,j:i\neq j} {v_iv_j^2}^\rcirclearrowleft-3\sum_{i,j,k\atop i,j,k\ \text{are distinct}} v_iv_jv_k.
\end{align*}
\end{theorem}

\begin{proof} By Theorem~\ref{thm03},
\begin{align*}
&P=a_1\sum_R S_R^3+a_2\sum_{R,i} S_R^2v_i^\rcirclearrowleft+a_3\sum_{R,i}S_R{v_i^2}^\rcirclearrowleft+a_4\sum_{i,j,R:ijR} S_Rv_iv_j^\rcirclearrowleft
\\
&+a_5\sum_{i,j,R:iRj} S_Rv_iv_j^\rcirclearrowleft+a_6\sum_{i} v_i^3+a_7\sum_{i,j:i\neq j} {v_iv_j^2}^\rcirclearrowleft+a_8\sum_{i,j,k\atop i,j,k\ \text{are distinct}} v_iv_jv_k
\end{align*}
for some non-zero $(a_1,\dots,a_8)\in \K^8$. The derivatives of $P$, being relations of $A^!$ must be orthogonal to quadratic relations of $A$ listed in Theorem~\ref{thm02}. Note that they are automatically orthogonal to all commutators as well as to products of distinct $S_R$. Orthogonality of $\frac{\partial P}{\partial S_R}$ to $S_{R'}t_{ij}$, $iR'j$ with $S_R\neq S_{R'}$ is also obvious.

Orthogonality of $\frac{\partial P}{\partial S_R}$ to $S_{R}t_{ij}$ with $iRj$ is easily seen to be equivalent to $4a_1+4a_2=0$, which yields
$$
a_2=-a_1.
$$

Orthogonality of $\frac{\partial P}{\partial S_R}$ to $t_{ij}t_{ik}$ with $ijR$ and $ikR$ reads (taking into account $a_2=-a_1$):
$$
84a_1+11a_3+57a_4-42a_5=0.
$$

Orthogonality of $\frac{\partial P}{\partial S_R}$ to $t_{ij}t_{ik}$ with $ijR$, $iRj$, or $iRj$, $ikR$, or $iRj$, $iRk$ reads (taking into account $a_2=-a_1$):
$$
-16a_1+11a_3-3a_4+8a_5=0.
$$

Orthogonality of $\frac{\partial P}{\partial v_i}$ to $S_Rt_{ij}$ with $iRj$ reads (taking into account $a_2=-a_1$):
$$
-4a_1+4a_3-2a_4+2a_5=0.
$$

Orthogonality of $\frac{\partial P}{\partial v_i}$ to $S_Rt_{jk}$ with $jRk$ and $i,j,k$ pairwise distinct reads (taking into account $a_2=-a_1$):
$$
-4a_1-a_3+3a_4+2a_5=0.
$$

The equations we have so far form a closed system of linear equations on $a_1,\dots,a_5$. The space of its solutions is one-dimensional and consists of $(a_1,\dots,a_5)$ satisfying
\begin{equation}\label{tttt}
\text{$a_2=-a_1$, $a_5=2a_1$, $a_3=a_4=0$.}
\end{equation}

We go on. Orthogonality of $\frac{\partial P}{\partial v_i}$ to $t_{ij}t_{ik}$ with $i,j,k$ pairwise distinct reads (taking into account (\ref{tttt})):
$$
-60a_1+16a_6-5a_7+5a_8=0.
$$

Orthogonality of $\frac{\partial P}{\partial v_i}$ to $t_{ij}t_{jk}$ with $i,j,k$ pairwise distinct reads (taking into account (\ref{tttt})):
$$
180a_1-4a_6+35a_7-15a_8=0.
$$

Finally, orthogonality of $\frac{\partial P}{\partial v_i}$ to $t_{jk}t_{jm}$ with $i,j,k,m$ pairwise distinct reads (taking into account (\ref{tttt})):
$$
20a_1+a_6+15a_8=0.
$$

Note that by now we have covered all options. The system of last three equations has one-dimensional space of solutions consisting of $(a_1,a_6,a_7,a_8)$ satisfying
\begin{equation}\label{tttt1}
\text{$a_6=\frac52a_1$, $a_7=-\frac{11}2a_1$, $a_8=-\frac32a_1$..}
\end{equation}

It remains to plug in some non-zero number as $a_1$ (we choose $a_1=2$ to kill denominators) and use (\ref{tttt}) and (\ref{tttt1}) to obtain the required formula for $P$.
\end{proof}

Now we consider the potential complex for $A^!$. It is easily seen to coincide with the Koszul complex. In order to prove that $A^!$ is Koszul (which is equivalent to Koszulity of $A$), it is enough to show that the potential complex for $A^!$  is exact. We do this by calculating the Hilbert series of $A^!$ and applying the following lemma.


The main part is the calculation of the Hilbert series. We will establish numerical Koszulity for $A^!$.

\begin{theorem}\label{theorem05}
 The Hilbert series of dual algebra $A^!$ for $A$ is
$$
H_{A^!}=\frac1{1-16t+16t^2-t^3}.
$$
\end{theorem}

First we find minimal series in the variety $P(16,3)$ of potential algebras with $16 $ generators and potential of degree $3$.

We use

\begin{proposition}\label{Prop1}
 The potential complex is exact in general position in P(n,m) if in P(n,m) there exists an exact potential complex.
 \end{proposition}

This could be shown in the same way as in \cite{IA}, Theorem 3.2.

Then we construct an exact potential(=Koszul)  complex.

{\bf Example} Consider a potential algebra $C=\K\langle x,y_1,\dots,y_s\rangle/\frac{\partial P}{\partial x},\frac{\partial P}{\partial y_j}$, $s\geq 2$ (we need $s=15$), with the potential

$$
P=xy_1y_2^\rcirclearrowleft+xy_2y_1^\rcirclearrowleft+\sum_{j=3}^s {xy_j^2}^\rcirclearrowleft.
$$

The relations are as follows:

\begin{align*}
&\frac{\partial P}{\partial x}=y_1y_2+y_2y_1+y_3^2+{\dots}+y_s^2;
\\
&\frac{\partial P}{\partial y_1}=xy_2+y_2x;
\\
&\frac{\partial P}{\partial y_2}=xy_1+y_1x;
\\
&\frac{\partial P}{\partial y_j}=xy_j+y_jx\ \ \text{for}\ j\geq3.
\end{align*}

It is easy to see that the only ambiguity $xy_1y_2$ resolves, so the relations form a quadratic Gr\"obner basis. Thus algebra $C$ is PBW and due to Priddy criterion \cite{Pr, Manin} it is Koszul.

Thus the exact potential complex according to proposition and the above Example is in general position in $P(n,m)$. Since generic series is minimal (standard arguments, see for example \cite{PP}), and exactness of potential complex defines the series uniquely via the recurrence relation, we know that minimal series in $P(16,3)$ is the same as in the Example and is equal to $\frac1{1-16t+16t^2-t^3}$.

So we have

\begin{proposition}\label{Prop2}
$H_{A^!}\geq \frac1{1-16t+16t^2-t^3}$.
\end{proposition}

Now we need to show that an equality in the opposite direction holds. For this we consider algebra $A^!$ over the field of characteristic 7. After coefficients of relations considered as elements of $\Z_7$ the pattern of Gr\"obner basis of relations becomes much simpler, namely the highest words can be made exactly the same as those in Example 1. And it is possible to check that these words form a Gr\"obner basis, hence the Hilbert series of $A^!$ is equal to
$$
H_{A^!_{\Z_7}}=\frac1{1-16t+16t^2-t^3}.
$$
But again standard arguments (as in \cite{PP}) show that after passing from $\Z$ to $\Z_p$, the series can only become bigger. Hence we have

\begin{proposition}\label{Prop3}
 $H_{A^!}\leq H_{A^!_{\Z_7}}=\frac1{1-16t+16t^2-t^3}$.
 \end{proposition}

Combining Proposition\ref{Prop2} and Proposition\ref{Prop3} we have that the Hilbert series of $A^!$ is
$$
H_{A^!}=\frac1{1-16t+16t^2-t^3}.
$$

To prove Theorem\ref{theorem00} we combine Theorem\ref{theorem05}, lemma\ref{NK} and note that the potential complex is exact in all but one terms.

\section{Acknowledgements }

 \medskip


We are grateful to IHES and MPIM for hospitality, support, and excellent research atmosphere. I would like to thank Yu.I.Manin for asking  this question and constant encouragement.
This work is funded by the ERC grant 320974, and partially supported by the project PUT9038.



  \normalsize

\vskip1truecm

\scshape

\noindent  Natalia Iyudu\

\noindent School of Mathematics

\noindent  The University of Edinburgh

\noindent James Clerk Maxwell Building

\noindent The King's Buildings

\noindent Peter Guthrie Tait Road

\noindent Edinburgh

\noindent Scotland EH9 3FD

\noindent E-mail addresses: \qquad {\tt n.iyudu@ed.ac.uk, \qquad{\tt n.joudu@yahoo.de}}

\vskip1truecm

\scshape

\end{document}